\input amstex
\input epsf.sty
\magnification=\magstep1
\documentstyle{amsppt}
\pagewidth{16 true cm} \pageheight{23.5 true cm}
\TagsOnRight
\NoRunningHeads

\document

\topmatter
\title The boundary of Young graph \\ with Jack edge multiplicities
\endtitle
\author
Sergei~Kerov, Andrei~Okounkov, and Grigori~Olshanski
\endauthor 
\thanks
The authors were partially supported by the grants 96-01-00676
(S.K.) and 95-01-00814 (A.O. and G.O.) of Russian Foundation for 
Basic Research.
\endthanks
\abstract
Consider the lattice of all Young diagrams ordered by inclusion,
and denote by $\Bbb{Y}$ its Hasse graph. Using the Pieri formula
for Jack symmetric polynomials, we endow the edges of the graph
$\Bbb{Y}$ with formal multiplicities depending on a real
parameter $\theta$. The multiplicities determine a potential
theory on the graph $\Bbb{Y}$. Our main result identifies the
corresponding Martin boundary with an infinite-dimensional simplex,
the ``geometric boundary'' of the Young graph $\Bbb{Y}$,
and provides a canonical integral representation for non-negative
harmonic functions.

For three particular values of the parameter, the theorem
specializes to known results: the Thoma theorem describing
characters of the infinite symmetric group, 
the Kingman's classification of partition structures,
and the description of spherical functions of the
infinite hyperoctahedral Gelfand pair.  
\endabstract
\endtopmatter

\subhead 1. The main result \endsubhead
We use \cite{13} as a basic reference for the standard notations
and terminology related to integer partitions and symmetric
functions (see also Appendix for basic background definitions).

Let $\Lambda$ denote the $\Bbb{R}$-algebra of symmetric functions
in infinitely many variables $x=(x_1,x_2,\dots)$. We denote by
$P_\mu(x;\theta)$ the Jack symmetric polynomial with the
parameter $\theta$, indexed by an integer partition $\mu$ (see
\cite{13}, VI.10, or Section A.4 for the definition; note that our
parameter $\theta$ is inverse to Macdonald's $\alpha=1/\theta$).
The family $\{P_\mu\}$ forms a linear basis in the symmetric
function algebra $\Lambda$. We aim to describe the set
$X(\theta)$ of algebra homomorphisms $\varphi:\Lambda\to\Bbb{R}$,
such that
$$
\varphi\big(P_\mu(\,\cdot\,;\theta)\big) \ge 0
\tag 1.1
$$
for a fixed value of $\theta$ and all integer partitions $\mu$.

Recall that the power sum symmetric polynomials
$p_m(x)=\sum_{j\ge1}x_j^m$ are algebraically independent and
generate the algebra $\Lambda$. As a result, every homomorphism
$\varphi:\Lambda\to\Bbb{R}$ is entirely determined by an
arbitrary real sequence $\varphi_1,\varphi_2,\ldots$ of its
values $\varphi_m=\varphi(p_m)$ at the polynomials $p_m$. The
main result of the present paper can be stated as follows.

\proclaim{Theorem A}
Assume that $\theta>0$, and let $\varphi:\Lambda\to\Bbb{R}$ be an
algebra homomorphism, such that $\varphi(p_1)=1$. Then the
following two conditions are equivalent:

(1)\; $\varphi(P_\mu(\,\cdot\,;\theta))\ge0$ for all integer
partitions $\mu$;

(2)\;there exists a pair of non-negative weakly decreasing
sequences
$$
\alpha=(\alpha_1\ge\alpha_2\ge\dots\ge0),\quad
\beta=(\beta_1\ge\beta_2\ge\dots\ge0),
$$
subject to the condition $\sum\alpha_j+\sum\beta_j\le1$,  and
$$
\varphi(p_m) =
\sum_{j=1}^\infty \alpha_j^m + (-\theta)^{m-1}
\sum_{j=1}^\infty \beta_j^m
$$
for all $m=2,3,\ldots$.

The sequences $\alpha,\beta$ are uniquely
determined by the homomorphism $\varphi$.
\endproclaim

This statement confirms a particular case of the Conjecture of
Section 7.3 in \cite{6}. In the special case $\theta=1$ the
Theorem is equivalent to Schoenberg's conjecture on totally
positive sequences, proved in \cite{1}, \cite{3}, and to a
theorem by E.~Thoma \cite{23} describing the characters of the
infinite symmetric group.

The limiting case $\theta=0$, in a different form, was studied by
Kingman \cite{11}, see also \cite{5}. Yet another particular case
where the Theorem was already known to be true is that of
$\theta=1/2$. In this case it admits a representation
theoretical interpretation, see \cite{18}, \cite{21}.

Theorem A can be restated in terms of discrete potential
theory. In this form it provides a Poisson-type integral
representation of non-negative harmonic functions on the Young
graph, with respect to an edge multiplicity function depending on
a real parameter $\theta$, see Theorem B in Section 5.

In a more general setup one starts with a Pieri-type formula
$$
p_1 \cdot P_\lambda = \sum_\nu \varkappa(\lambda,\nu)\, P_\nu
$$
where $\{P_\lambda\}$ is a linear homogeneous basis in the
algebra $\Lambda$. The coefficients $\varkappa(\lambda,\nu)$
(assumed to be non-negative) determine formal multiplicities of
edges of the Young graph $\Bbb{Y}$. One looks for an integral
representation of non-negative functions $\varphi$ on the set
$\Bbb{Y}$, harmonic in the sense that
$$
\varphi(\lambda) =
\sum_\nu \varkappa(\lambda,\nu)\, \varphi(\nu).
$$

Theorem B solves the problem for the basis of Jack symmetric
polynomials $P_\lambda(x;\theta)$. The proof relies substantially
on a formula for $\theta$-dimension of skew Young diagrams
obtained in \cite{20}. We also use many ideas and constructions
of papers \cite{5-9}, \cite{24} though the present paper can be
read independently.

A similar integral representation is known for Schur
$Q$-functions (a special case of Hall -- Littlewood symmetric
polynomials,  see \cite{13}, III.8). The corresponding graph is
that of shifted Young diagrams (see Nazarov \cite{15} and Ivanov
\cite{4}), it describes the branching of projective
representations of symmetric groups.

A part of the conjecture in \cite{6} related to general
Macdonald symmetric polynomials still remains open.

The plan of the paper is as follows. In Section 2 we recall the
Martin boundary construction for graded graphs with
multiplicities. In Sections 3, 4 we present the well-known
integral representations of non-negative harmonic functions for
two special cases of Jack graphs: the Young lattice and the
Kingman's graph. The general case of arbitrary Jack graphs is
considered in Section 5. In Sections 6, 7 we derive two basic
ingredients in the proof of our main Theorems A and B. The proof
itself is given in Section 8. There is also an appendix where we
recall necessary background definitions related to integer
partitions and symmetric polynomials.

\subhead 2. The Martin boundary of a branching graph \endsubhead
In this paper we consider a particular case of the following
general problem: find the Martin boundary of an infinite graph
$\Gamma$ with respect to a given multiplicity function on the set
of its edges.

More precisely, the graph $\Gamma$ is graded 
$$
\Gamma=\bigcup_{n=0}^\infty\Gamma_n\,,
$$
 the vertices in $\Gamma_n$ representing the admissible
states of a process at the moment $n$. The kinematics of the
process is described by the edges of the graph pointing out from
a vertex $\lambda\in\Gamma_n$ to some vertex $\nu$ at the next
level $\Gamma_{n+1}$. We write $\lambda\nearrow\nu$ to show that
the pair $(\lambda,\nu)$ is an oriented edge of the graph
$\Gamma$.  The process always starts with a distinguished vertex
$\varnothing\in\Gamma_0$ at the zero level of the graph $\Gamma$.

Another piece of information for our general problem is provided
by a {\it multiplicity function} $\varkappa$ on the set of edges
of $\Gamma$, with real positive values. A pair
$(\Gamma,\varkappa)$ consisting of a graded graph $\Gamma$ and a
multiplicity function $\varkappa$ will be referred to as a {\it
branching graph} (the origin of the term will be motivated by an
example in Section 3).

A function $\varphi:\Gamma\to\Bbb{R}$ defined on the set of
vertices in $\Gamma$ is called {\it harmonic}, if the following
variant of the ``mean value theorem'' holds for all vertices
$\lambda\in\Gamma$:
$$
\varphi(\lambda) = \sum_{\nu:\lambda\nearrow\nu}
\varkappa(\lambda,\nu)\; \varphi(\nu).
\tag 2.1
$$
We are interested in the space $\Cal{H}$ of all non-negative
harmonic functions normalized at the vertex $\varnothing$ by the
condition $\varphi(\varnothing)=1$. The topology in the space
$\Cal{H}$ is that of pointwise convergence. It is clear that
$\Cal{H}$ is convex and compact.

A general solution for such a problem is based on the Martin
boundary construction (see, for instance, \cite{2}). One starts
with the {\it dimension function} $\dim(\mu,\nu)$ defined
recurrently by the formulae $\dim(\mu,\mu)\equiv1$,
$$
\dim(\mu,\nu) = \sum_{\lambda:\lambda\nearrow\nu}
\dim(\mu,\lambda)\, \varkappa(\lambda,\nu),
\tag 2.2
$$
and $\dim(\mu,\nu)=0$, if there is no oriented path $t$ from $\mu$
to $\nu$. Associating a weight
$w(t)=\prod\varkappa(\lambda_{i-1},\lambda_i)$ with every such
path $t=(\lambda_0=\mu,\lambda_1,\,\ldots,\lambda_m=\nu)$, one
can write the dimension function as a sum of weights over all
oriented paths between the vertices $\mu$ and $\nu$, that is,
$\dim(\mu,\nu)=\sum w(t)$.

From the point of view of potential theory,
$G(\mu,\nu)=\dim(\mu,\nu)$ is the Green function with respect to
``Laplace operator''
$$
(\Delta\varphi)(\lambda) = -\varphi(\lambda) +
\sum_{\nu:\lambda\nearrow\nu}
\varkappa(\lambda,\nu)\,\varphi(\nu).
$$
This means that if $\varphi_\nu(\mu)=G(\mu,\nu)$ for a fixed vertex
$\nu$, then $-(\Delta\varphi_\nu)(\mu)=\delta_{\mu\nu}$ for all
$\mu\in\Gamma$. The ratio
$$
K(\mu,\nu) = \frac{\dim(\mu,\nu)}{\dim\nu}
\tag 2.3
$$
is usually called the {\it Martin kernel}.

Consider the space $\Cal{F}$ of all functions
$f:\Gamma\to\Bbb{R}$ with the topology of pointwise convergence,
and let $\widetilde{E}$ be the closure of the subset
$\widetilde{\Gamma}\subset\Cal{F}$ of functions
$\mu\mapsto{K}(\mu,\nu)$, $\nu\in\Gamma$. Since those functions
are uniformly bounded, $0\le{K}(\mu,\nu)\le1$, the space
$\widetilde{E}$ (called the {\it Martin compactification}) is
indeed compact. One can easily check that
$\widetilde{\Gamma}\subset\widetilde{E}$ is a dense open subset
of $\widetilde{E}$. Its boundary
$E=\widetilde{E}\setminus\widetilde{\Gamma}$ is called the {\it
Martin boundary} of the branching graph $(\Gamma,\varkappa)$.

By definition, the Martin kernel (2.3) may be extended by
continuity to the function
$K:\Gamma\times\widetilde{E}\to\Bbb{R}$, and we keep denoting it
by the same letter $K$. For each boundary point $\omega\in{E}$
the function $\varphi_\omega(\mu)=K(\mu,\omega)$ is non-negative,
harmonic, and normalized. Moreover, all such functions have
an integral representation similar to the classical Poisson
integral for non-negative harmonic functions in the disk.

\proclaim{Theorem \rm (cf. \cite{2})}
Every normalized non-negative harmonic function
$\varphi\in\Cal{H}$ admits an integral representation
$$
\varphi(\lambda) = \int_E K(\lambda,\omega)\;M(d\omega).
\tag 2.4
$$
where $M$ is a probability measure. For every probability measure
$M$ on $E$ the integral (2.4) provides a non-negative harmonic
function $\varphi\in\Cal{H}$.
\endproclaim

All {\it indecomposable} (i.e., {\it extreme}) functions in
$\Cal{H}$ can be represented in the form
$\varphi_\omega(\mu)=K(\mu,\omega)$, for appropriate boundary
point $\omega\in{E}$, and we denote by $E_{min}$ the
corresponding subset of the boundary $E$. It is known that
$E_{min}$ is a non-empty $G_\delta$ subset of $E$. One can always
choose the measure $M$ in the integral representation (2.4) to be
supported by $E_{min}$. Under this assumption, the measure $M$
representing a function $\varphi\in\Cal{H}$ via (2.4) is unique.

In all the examples considered in the present paper the linear
span of the functions $\omega\mapsto{K}(\mu,\omega)$,
$\mu\in\Gamma$ is dense in the space $C(E)$ of continuous
functions on the boundary $E$. As a consequence, the minimal
boundary $E_{min}$ will always coincide with the entire Martin
boundary $E$.

Given a concrete example of a branching graph, one looks for an
appropriate ``geometric'' description of the abstract Martin
boundary. For instance, the Martin boundary of a domain in the
complex plane interior to a Jordan curve can be identified with
this curve, its geometric boundary. For a general simply
connected domain the Martin boundary is more complicated; its
geometric version is provided by the Perron -- Carath\'eodory
theory of ``boundary elements'', see \cite{14}, Chapter V, \S3.

The purpose of the present paper is to give an explicit integral
representation of the form (2.4) for the Young graph $\Bbb{Y}$
endowed with some edge multiplicities
$\varkappa_\theta(\lambda,\nu)$ depending on a real parameter
$\theta$. The multiplicities arise in the Pieri formula for the
Jack symmetric polynomials $P_\mu(x;\theta)$. We refer to the
corresponding branching as to the Jack branching graph, and
denote it as $\Bbb{J}(\theta)$. Presently, we proceed with a
couple of well-studied concrete examples of such graphs. Both
examples will appear as particular instances in the family
$\Bbb{J}(\theta)$.

\subhead 3. Example: the Young lattice \endsubhead
By definition, the vertices of the {\it Young graph} $\Bbb{Y}$
are Young diagrams (representing partitions of natural numbers).
The edges join the pairs of diagrams which only differ by a
single box (i.e., if the second diagram in a pair $(\lambda,\nu)$
covers the first one with respect to inclusion order).

$${\epsffile{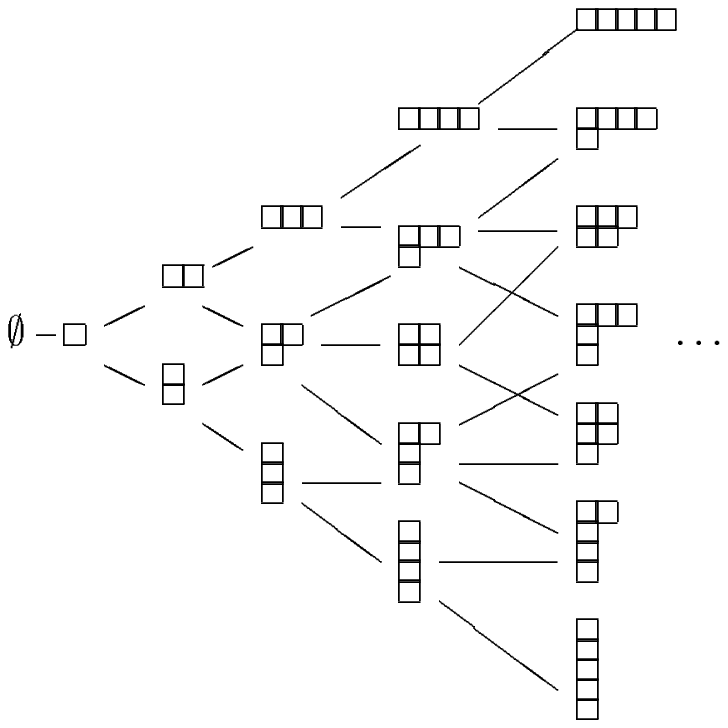}}$$

\centerline{\tenrm Fig.~1\quad The Young graph.}
\bigskip

The graph $\Bbb{Y}$ plays a fundamental role in representation
theory of finite and infinite symmetric groups, see Vershik and
Kerov \cite{7}, \cite{9}, \cite{24}. The reason is that the set
$\Bbb{Y}_n$ of Young diagrams with $n$ boxes labels characters of
irreducible representations of the symmetric group $\frak{S}_n$
of degree $n$, and the edges of the Young graph describe the
branching of characters of $\frak{S}_n$ when restricted to
the subgroup $\frak{S}_{n-1}$. The multiplicity function is
trivial: $\varkappa\equiv1$ for all edges. The first five levels
of the Young graph are represented in Fig.~1.

The set of all central positive definite functions on the
infinite symmetric group $\frak{S}_\infty=\bigcup\frak{S}_n$,
normalized at the identity, can be naturally identified with the
set $\Cal{H}$ of positive harmonic functions on the Young graph.
The {\it characters} of the group $\frak{S}_\infty$, i.e., the
indecomposable normalized central positive-definite functions,
correspond to the points of the set $E_{min}=E$.

The space of characters of the infinite symmetric group
$\frak{S}_\infty$ was found by E.~Thoma in the pioneering paper
\cite{23}. In equivalent terms of potential theory on the Young
graph his result can be stated as follows.

\smallskip\noindent{\bf Definition.}
The space $\Omega$ of all pairs $\omega=(\alpha,\beta)$ of
weakly decreasing sequences of non-negative real numbers,
$$
\alpha = (\alpha_1\ge\alpha_2\ge\ldots\ge\alpha_n\ge\ldots\ge0),
\quad
\beta = (\beta_1\ge\beta_2\ge\ldots\ge\beta_n\ge\ldots\ge0),
\tag 3.1
$$
such that 
$$
\gamma=1-\sum\alpha_n-\sum\beta_n\ge0\,,
$$
is called
the {\it Thoma simplex}. The topology of $\Omega$ is that of
pointwise convergence.

\smallskip
Given a Young diagram $\nu$, we denote by
$d=\max\{j:\nu_j\ge{j}\}$ the length of its diagonal, and by
$$
a_j = \nu_j-j,\qquad b_j = \nu_j'-j+1;\qquad j = 1,\,\ldots,d
\tag 3.2
$$
the so called {\it Frobenius parameters} of the diagram $\nu$.
We associate with $\nu$ a point
$$
\omega_\nu = \left(
\frac{a_1}n\;,\dots,\frac{a_d}n\;,0,\dots\;;\;
\frac{b_1}n\;,\dots,\frac{b_d}n\;,0,\dots
\right)
\tag 3.3
$$
of the Thoma simplex $\Omega$. Remark that a Young diagram $\nu$
is entirely determined by the point $\omega_\nu$, along with the
value of $n=|\nu|$. Thus, we may identify the set of non-empty
Young diagrams with the subset
$$
\widetilde{\Bbb Y} =
\bigcup_{n=1}^\infty \bigcup_{\nu\in\Bbb{Y}_n}
\Big({1 \over n},\;\omega_\nu\Big)
\tag 3.4
$$
of the product space $[0,1]\times\Omega$. The set
$\widetilde{\Bbb{Y}}$ is discrete (i.e., each of its points has a
neighborhood in $[0,1]\times\Omega$ free of other points of
$\widetilde{\Bbb{Y}}$), and its boundary $\{0\}\times\Omega$ is
homeomorphic to Thoma simplex. By this reason we shall call
$\Omega$ the {\it geometric boundary of the Young lattice}.

We are now in a position to describe the Martin boundary of the
Young graph.

\proclaim{Theorem \rm (cf. \cite{23}, \cite{24}, \cite{16})}
The Martin boundary $E$ of the Young graph is naturally
homeomorphic to its geometric boundary, the Thoma simplex
$\Omega$. Every non-negative normalized harmonic function
$\varphi\in\Cal{H}$ on the Young graph admits a unique
integral representation of the form
$$
\varphi(\lambda) = \int_\Omega s_\lambda(\omega)\; M(d\omega).
\tag 3.5
$$
The Martin kernel $K(\lambda,\omega)=s_\lambda(\omega)$ is given
by the extended Schur functions (see Section A.3 for the
definition of extended symmetric functions).
\endproclaim

The boundary $\Omega$ of the Young graph $\Bbb{Y}$ has quite a
number of interpretations. For instance, the equation (A.2)
establishes a bijection between $\Omega$ and the set of all
totally positive real sequences $\{h_n\}_{n=0}^\infty$. Such
sequences were studied  by I.~Schoenberg and his
school \cite{1}. Schoenberg's conjecture about the
structure of totally positive sequences 
was proved by A.~Edrei \cite{3} long before the paper \cite{23} was published, and
independently of the representation theory of the group
$\frak{S}_\infty$. 

Yet another interpretation is based on the observation that the
formula $$ \psi(h_n) = h_n(\alpha,\beta);\qquad (\alpha;\beta)
\in \Omega, \quad n=1,2,\ldots $$ determines a general
homomorphism $\psi:\Lambda\to\Bbb{R}$ of the ring $\Lambda$ of
symmetric polynomials, which is non-negative on the basis of
Schur functions, $\psi(s_\lambda)\ge0$ for all
$\lambda\in\Bbb{Y}$, see \cite{23}.

Young graph governs the branching of irreducible characters of
the symmetric groups $\frak{S}_n$. There exists a similar
branching graph $\Bbb{K}$ responsible for the branching of
conjugacy classes in these groups.

The underlying graph of $\Bbb{K}$ coincides with the Young graph
$\Bbb{Y}$, but this time the multiplicities of edges are
non-trivial. More precisely,
$$
\varkappa(\lambda,\nu) = m_k(\nu),
\tag 4.1
$$
where $k$ denotes the length of the row of the diagram $\nu$
containing a box to be removed from $\nu$ in order to obtain
$\lambda$, and $m_k(\nu)$ is the total number of rows of length
$k$ in the diagram $\nu$ (see Fig.~2).

Motivated by a problem of population genetics, J.~F.~C.~Kingman
introduced in \cite{11} the notion of {\it partition structure}. It
can be defined as a family $\{M_n\}_{n=1}^\infty$ where $M_n$ is
a central (i.e., constant on conjugacy classes) probability
distribution on the symmetric group $\frak{S}_n$. The family is
assumed to be {\it coherent} in the sense that $d(M_{n+1})=M_n$,
for all $n=1,2,\ldots$, where the permutation $d(w)\in\frak{S}_n$
is obtained from $w\in\frak{S}_{n+1}$ by crossing the element
$n+1$ out of its cycle in $w$.

The description of Martin boundary for the branching graph
$\Bbb{K}$ is equivalent to the Kingman's classification of
partition structures, and can be stated as follows.

\proclaim{Theorem \rm (cf. \cite{11} and \cite{5})}
Every non-negative function $\varphi$ on the set of Young
diagrams, harmonic with respect to the multiplicities (4.1),
admits a unique integral representation
$$
\varphi(\lambda) = \int_E m_\lambda(\alpha)\; M(d\alpha),
\tag 4.2
$$
where $M$ is a probability measure, and the boundary $E=E_{min}$
is that part of the Thoma simplex (3.2) for which $\beta_n=0$ for
all $n=1,2,\ldots$. The Martin kernel $m_\lambda(\alpha)$ is
given by the extended\footnote{This tiny extension concerns only
the $\gamma$ parameter; all $\beta$-parameters are supposed to
vanish, see Section A.3.} monomial symmetric functions.
\endproclaim

$${\epsffile{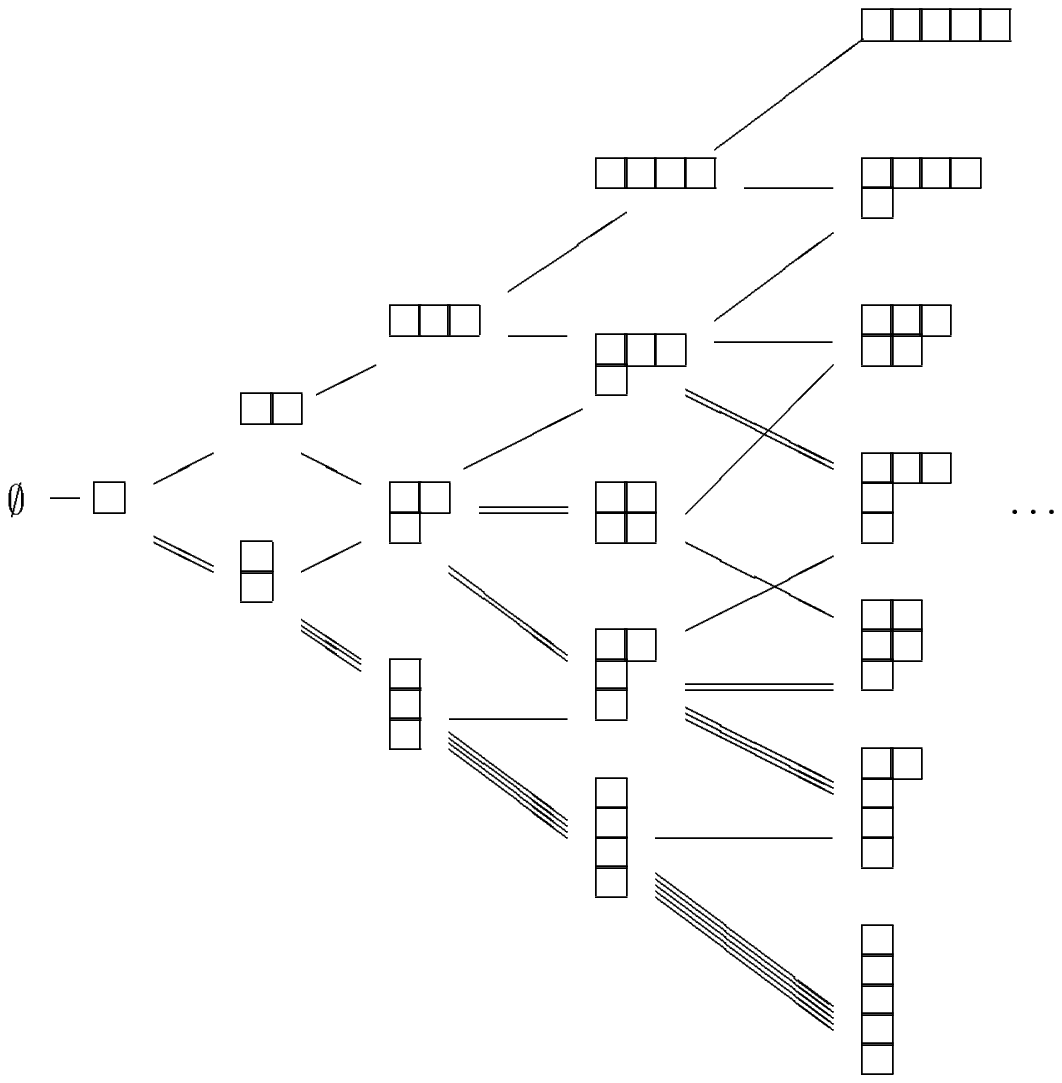}}$$

\centerline{\tenrm Fig.~2\quad The Kingman's branching graph.}
\bigskip

Kingman's result is much simpler than Thoma theorem. See
\cite{5} for another proof of this result.

\subhead 5. Jack's graphs and their boundaries \endsubhead
In this Section we introduce a one-parameter family of
multiplicity functions on the Young graph, interpolating between
the trivial multiplicities and those of Kingman. Our main Theorem
B identifies the corresponding boundary and the Martin kernel.

We refer to \cite{13} for the theory of Jack symmetric polynomials
(see also equations (A.5) -- (A.7) for some definitions). We shall
use the simplest particular case of Pieri formula for Jack
symmetric polynomials which reads
$$
p_1(x)\, P_\lambda(x;\theta) = \sum_{\nu:\lambda\nearrow\nu}
\varkappa_\theta(\lambda,\nu)\, P_\nu(x;\theta),
\tag 5.1
$$
where the multiplicities $\varkappa_\theta(\lambda,\nu)$ are
given by an explicit formula
$$
\varkappa_\theta(\lambda,\nu) = \prod_b \frac
{\big(a(b) + (l(b)+2)\theta\big)\big(a(b) + 1 + l(b)\theta\big)}
{\big(a(b)+1+(l(b)+1)\theta\big)\big(a(b)+(l(b)+1)\theta\big)}\;.
\tag 5.2
$$
Here $b$ runs over all boxes in the $j$-th column of the diagram
$\lambda$, provided that the new box $b_0=\nu\setminus\lambda$
belongs to the $j$-th column of $\nu$. The number
$a(i,j)=\lambda_i-j$ is called the {\it arm length}, and
$l(i,j)=\lambda_j'-i$ is the {\it leg length} of a box $b=(i,j)$
in the diagram $\lambda$. (See \cite{13}, Ch. VI, (10.10) and
(6.24.iv).)

It can be checked that the multiplicities (5.2) tend, as
$\theta\to0$, to those of Kingman (4.1) (the limit behavior of
(5.2) only depends on the boxes $b$ with $a(b)=0$).

\proclaim{Lemma 5.1} The coefficients (5.2) are all positive iff
$\theta\ge0$. \endproclaim

\demo{Proof}
Clearly, the condition is sufficient. Since
$\varkappa_\theta\big((1^n),(1^{n+1})\big)=(1+n)/(1+n\theta)$, it
is necessary, too.
\qed\enddemo

\smallskip\noindent {\bf Definition.}
We denote by $\Bbb{J}(\theta)=(\Bbb{Y},\varkappa_\theta)$ the
Young graph $\Bbb{Y}$ with the edge multiplicities (5.2). By
abuse of terminology we call the branching diagram
$\Bbb{J}(\theta)$ the {\it Jack graph} with the parameter
$\theta$. We shall always assume that $\theta\ge0$, so that the
multiplicities are positive.

\smallskip
The initial part of the graph $\Bbb{J}(\theta)$ is represented in
Fig.~3. Note that $\Bbb{J}(0)$ coincides with the Kingman's
branching of Section 4, and that $\Bbb{J}(1)$ is the ordinary
Young graph with trivial multiplicities.

Our main concern in this paper is to describe all non-negative
harmonic functions on the graph $\Bbb{J}(\theta)$.

Assuming that $\theta>0$, we shall prove that the Martin
boundary of the graph $\Bbb{J}(\theta)$ does not depend on the
the parameter $\theta$ and may be naturally identified with the
Thoma simplex $\Omega$ -- the geometric boundary of the Young
graph. On the contrary, the Martin kernel does depend on
$\theta$. In order to describe this kernel, we associate to each
point $(\alpha;\beta)\in\Omega$ an algebra homomorphism
$\varphi_{\alpha,\beta}:\Lambda\to\Bbb{R}$ defined by the
formulae $\varphi_{\alpha,\beta}(p_1)\equiv1$ and
$$
\varphi_{\alpha,\beta}(p_m) = \sum_{j=1}^\infty \alpha_j^m +
(-\theta)^{m-1}\, \sum_{j=1}^\infty \beta_j^m;
\qquad m=2,3,\ldots\;.
\tag 5.3
$$

$${\epsffile{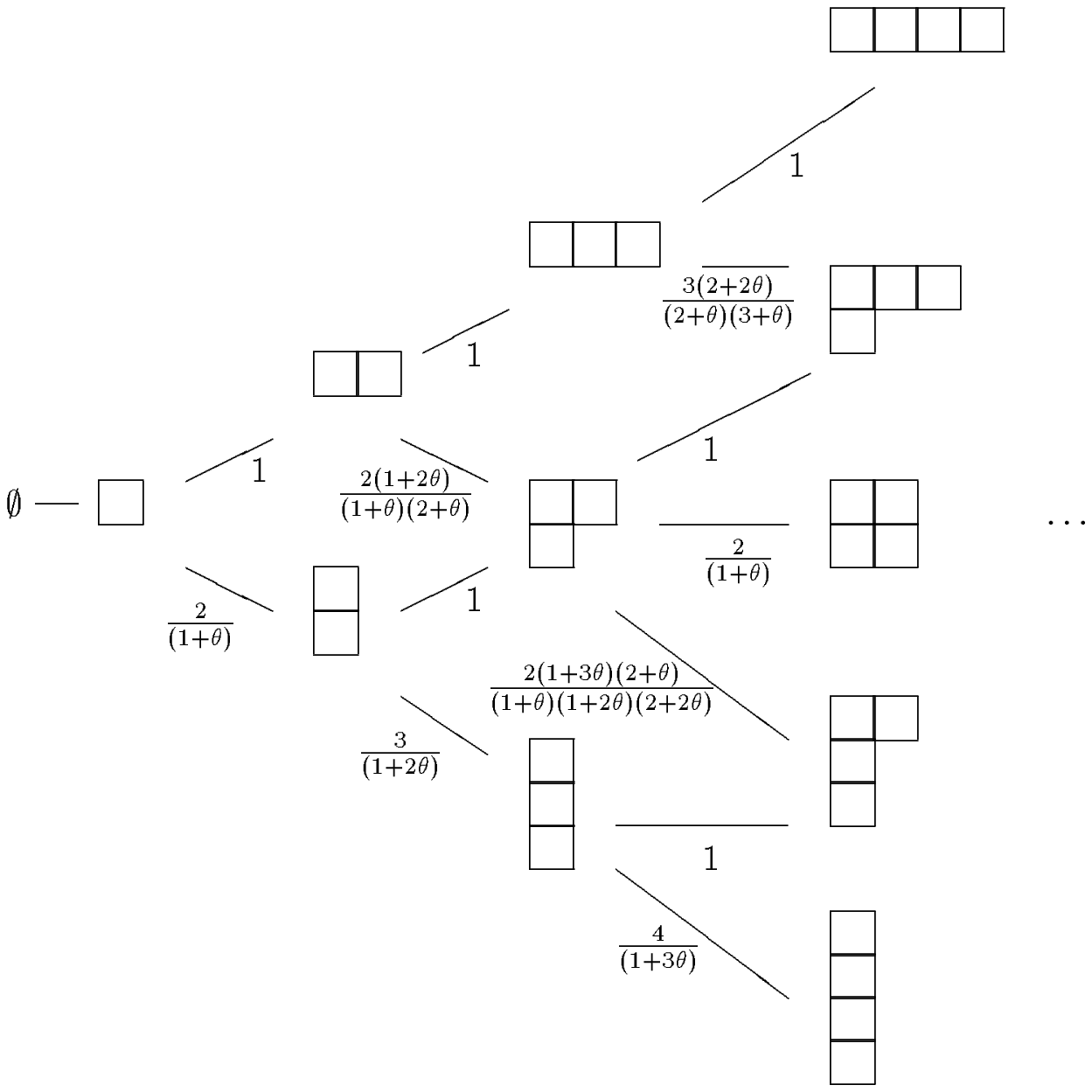}}$$

\centerline{\tenrm Fig.~3\quad The Jack's branching graph.}

\smallskip\noindent {\bf Definition.}
Given a symmetric polynomial $f\in\Lambda$, we consider a
function
$$
f(\alpha;\beta|\theta) = \varphi_{\alpha,\beta}(f),\qquad
(\alpha;\beta)\in\Omega,
\tag 5.4
$$
defined on the Thoma simplex $\Omega$, and we call it {\it
$\theta$-extension} of the polynomial $f$.

\smallskip
Note that if $\sum\alpha_i=1$ then all $\beta$-variables
vanish, $\beta=\bold0$, and the function
$f(\alpha)=f(\alpha;\bold{0})$ is identical with the naive
evaluation of the polynomial $f$ at the sequence
$\alpha=(\alpha_1,\alpha_2,\ldots)$. In case of $\theta=1$ the
$\theta$-extended functions coincide with the extended functions
of Section A.3.

\proclaim{Theorem B}
Assume that $\theta>0$, and let $\omega=(\alpha;\beta)$
denote a point of Thoma simplex $\Omega$. Then the
$\theta$-extension
$$
K(\mu;\omega) = \varphi_{\omega}\big(P_\mu(\,\cdot\,;\theta)\big),
\tag 5.5
$$
of the Jack symmetric polynomial $P_\mu(x;\theta)$ provides the
Martin kernel of the Jack graph $\Bbb{J}(\theta)$. The integral
representation
$$
\varphi(\mu) = \int_\Omega K(\mu;\omega)\, M(d\omega)
\tag 5.6
$$
establishes a one - to - one correspondence between the space
$\Cal{H}$ of normalized, non - negative functions
$\varphi:\Bbb{Y}\to\Bbb{R}$, harmonic with respect to Jack
multiplicity function $\varkappa_\theta$, and the space $\Cal{M}$
of probability measures $M$ on the Thoma simplex $\Omega$.
\endproclaim

This Theorem will be proved in Section 8 as a corollary of some
preliminary work in Sections 6, 7.

Besides the particular cases $\theta=0,1$ considered in Sections
3, 4 there is an interesting special case $\theta=1/2$. The
elements of $X(1/2)$ correspond to spherical unitary
representations of the Gelfand pair consisting of the ``even''
infinite symmetric group
$\frak{S}(2\infty)=\varinjlim\frak{S}_{2n}$ and its
hyperoctahedral subgroup $\varinjlim
\frak{S}(n)\ltimes\Bbb{Z}_2^n$. The
description of the Martin boundary $X(1/2)$ was obtained
(in representation-theoretical terms) in \cite{18}, see also
\cite{21}.

We complete the present Section with a couple of properties of
$\theta$-extended symmetric functions to be used later on.

\proclaim{Lemma 5.2}
Let $f\in\Lambda$ be a symmetric polynomial. Then its
$\theta$-extension $\omega\mapsto{f}(\omega|\theta)$ is
continuous on the Thoma simplex $\Omega$.
\endproclaim

\demo{Proof}
Trivially, the function $p_1(\alpha;\beta|\theta)\equiv1$ is
continuous. Since a sequence $\alpha$ is decreasing, we obtain
inequalities $k\alpha_k\le{\alpha}_1+\ldots+\alpha_k\le1$ and
hence $\alpha_k\le{1/k}$. A similar inequality $\beta_k\le1/k$
holds for a sequence $\beta$. Given a positive $\varepsilon>0$
and $m=2,3,\ldots$, there exists such an $N$ that
$$
\Big|\sum_{k=N}^\infty \alpha_k^m\Big| +
|\theta|^{m-1}\, \Big|\sum_{k=N}^\infty \beta_k^m\Big|
\le \varepsilon,
$$
i.e., the $N$-th tail of the series for the functions
$p_m(\alpha;\beta|\theta)$ is uniformly small on $\Omega$. 
It follows that the functions $p_m(\alpha;\beta|\theta)$ are continuous. All
other $\theta$-extended symmetric functions are polynomials in
the generators $p_1,p_2,\ldots$, and the Lemma follows.
\qed\enddemo

\proclaim{Lemma 5.3}
Assume that $\theta\ne0$. Then the algebra of $\theta$-extended
functions $\omega\mapsto{f}(\omega|\theta)$, $f\in\Lambda$, is
uniformly dense in the space $C(\Omega)$ of continuous functions
on the Thoma simplex.
\endproclaim

\demo{Proof}
By the Stone -- Weierstrass theorem, it suffices to show that the
functions $\omega\mapsto{p}_m(\omega|\theta)$ separate the points
of the Thoma simplex $\Omega$. One easily derives from the
definition (5.3) the generating series
$$
\sum_{m=1}^\infty {p_m(\alpha;\beta|\theta) \over z^m} =
\sum_{k=1}^\infty {\alpha_k \over z - \alpha_k} +
\sum_{k=1}^\infty {\beta_k \over z + \theta\beta_k} +
{\gamma \over z}
$$
for $\theta$-extended power sum polynomials. Hence, the values of
$\alpha_k$, $-\theta\beta_k$ can be restored as non-zero poles of
this series.
\qed\enddemo

\def\li #1 #2 #3 {\vrule height #1 pt width #2 pt depth #3 pt}

\def\hl #1 #2 #3 #4 {\rlap{ \kern #1 pt \raise -#2 pt \hbox{ \li #4 #3 0 }}}

\def\vl #1 #2 #3 #4 {\rlap{ \kern #1 pt \raise -#2 pt \hbox{ \li #3 #4 0 }}}

\def\wr #1 #2 #3 {\rlap{ \kern #1 pt \raise -#2 pt \hbox{\tenrm #3 }}}

\def\picture #1 #2 #3 \endpicture
   { \bigskip \centerline{\hbox to #1 pt
   {\nullfont #3 \hss}} \nobreak \bigskip \centerline{#2} \bigskip}

\def\th{\theta}
\def\kk{\kappa_\th}
\def\l{\lambda}
\def\k{\varkappa_\th}
\def\tht{\thetag}

\remark{Remark} 
Recall that the polynomial 
$$
J_\mu(x;\th)= H'(\mu) P_\mu(x;\th)
$$
where $H'(\mu)=\prod_{s\in\mu}(a(s)+\th l(s) + \th)$
is called the {\it integral form} of the Jack polynomial
\cite{13}. Define the coefficients $\kk(\mu,\l)$ as
those in the Pieri formula for polynomials 
$J_\mu$, 
$$
p_1 \cdot J_\mu = \sum_\l \kk(\mu,\l)\, J_\l\,.
$$
Since
$$
\kk(\mu,\l)= \frac{H'(\mu)}{H'(\l)} \k(\mu,\l), 
\tag 5.7
$$
the coefficients lead to an equivalent Martin
boundary problem for the Young graph. On the other
hand, the coefficients $\kk(\mu,\l)$ have the
following especially symmetric expression. 

\picture 100 {Fig.~4\quad Inner and outer corners of a diagram.}
\hl 0 0 100 1
\vl 0 70 70 1
\hl 0 60 18.5 0.5
\vl 20 58 8 0.5
\hl 20 50 8.5 0.5 
\vl 30 48 28 0.5
\hl 30 20 38.5 0.5 
\vl 70 18 8 0.5
\hl 70 10 18.5 0.5 
\vl 90 8 8 0.5
\hl 30 30 10 0.5 
\vl 40 30 10 0.5
\wr 26 22 {$\bigstar$}
\wr 18 52 {$\bullet$}
\wr -2 62 {$\bullet$}
\wr 68 12 {$\bullet$}
\wr 88 2 {$\bullet$}
\wr 28 52 {$\circ$}
\wr 18 62  {$\circ$}
\wr 68 22 {$\circ$}
\wr 88 12 {$\circ$}
\wr -2 76 {$x$}
\wr 103 2 {$y$}
\wr 42 33 {$\lambda/\mu$}
\endpicture
Fig.~4 represents inner and outer corners of the diagram $\mu$. 
The corner of $\mu$ where 
the box $\lambda/\mu$ is attached to $\mu$ is 
marked by $\bigstar$.
One easily checks that the products in \tht{5,2}, \tht{5.7}
telescope to 
$$
\kk(\mu,\l)=\frac 1{\th} \, \frac
{\dsize \prod_{
\text{outer corners $\circ$ of $\mu$}}
\, r_\th(\circ,\bigstar)}
{\dsize \prod_{
\text{inner corners $\bullet$ of $\mu$}}
\, r_\th(\bullet,\bigstar)}\,, \tag 5.8 
$$
where the {\it $\th$-axial distance} $r_\th$ 
between two  points on
the $(x,y)$ plane is defined by
$$
r_\th\big((x_1,y_1)-(x_2,y_2)\big):=
(y_1-y_2)- \th(x_1-x_2)\,,
$$
and the product in the denominator of \tht{5.8}
ranges over all inner corners distinct from the
corner $\bigstar$.

For the coefficients $\k(\mu,\l)$ there is a similar but 
less symmetric expression. 
\endremark

\subhead 6. Shifted Jack polynomials and the dimension formula
\endsubhead
In this Section we establish the first major ingredient in our
proof of Theorem B -- the explicit formula for the Martin
kernel of the Jack graph $\Bbb{J}(\theta)$, found in the paper
\cite{20}. To this end we survey some necessary facts from
\cite{20} concerning $\theta$-shifted symmetric polynomials in
general, and $\theta$-shifted Jack symmetric polynomials in
particular.

Denote by $\Lambda^\theta(n)$ the subalgebra of
$\Bbb{R}[x_1,\dots,x_n]$ formed by the polynomials sym\-me\-tric
in the `shifted' variables $x'_j=x_j-\theta j$, $j=1,\dots,n$.
We define the projection map
$\Lambda^\theta(n)\to\Lambda^\theta(n-1)$ as the specialization
$x_n=0$, and note that this projection preserves the filtration
defined by ordinary degree of polynomials. The projective limit
$\Lambda^\theta=\varprojlim\Lambda^\theta(n)$ in the category of
filtered algebras is referred to as the {\it algebra of
$\theta$-shifted symmetric polynomials}. The degree of an element
$F\in\Lambda^\theta$ is denoted as $\deg F$.

Each element $F\in\Lambda^\theta$ can be evaluated at any
sequence $x=(x_1,x_2,\dots)$ with finitely many non-zero terms. In
particular, one can evaluate $F$ at any integer partition
$\nu\in\Bbb{Y}$, which will be important in what follows.

For each polynomial $F\in\Lambda^\theta$, we denote by
$[F]\in\Lambda$ its {\it leading symmetric term} which is a
homogeneous symmetric polynomial of degree $\deg F$. The map
$F\mapsto[F]$ provides a canonical isomorphism of the graded
algebra $\operatorname{gr}\Lambda^\theta$ associated to the
filtered algebra $\Lambda^\theta$, onto the symmetric function
algebra $\Lambda$. Assuming that the leading terms
$[F_1],[F_2],\ldots$ of a sequence
$F_1,F_2,\ldots\in\Lambda^\theta$ generate the algebra $\Lambda$,
one readily derives that the latter polynomials generate the
algebra $\Lambda^\theta$.

For an elementary example of a $\theta$-shifted symmetric
function, consider a polynomial
$$
p^*_m(x;\theta) = \sum_{j=1}^\infty
\Big((x_j - \theta j)^m - (-\theta j)^m\Big),
\tag 6.1
$$
a shifted analog of the power sum symmetric function
$p_m=\sum{x}_j^m$. Since $[p_m^*]=p_m$, the polynomials
$p_1^*,p_2^*,\ldots$ generate the algebra $\Lambda^\theta$.

We are interested in $\theta$-shifted counterparts
$P_\mu^*(x;\theta)$ of the Jack symmetric polynomials
$P_\mu(x;\theta)$. These polynomials were studied
in \cite{12}, \cite{17} and \cite{20}; we refer to these papers
for a detailed exposition (see also \cite{20} which is dedicated
to the special case $\theta=1$).

The polynomial $P_\mu^*(x;\theta)$ can be characterized as the
unique element of the algebra $\Lambda^\theta$, such that $\deg
P_\mu^*=|\mu|$ and
$$
P_\mu^*(\lambda;\theta) = \cases
H(\mu), &\lambda=\mu\,,\\
0, &\mu\not\subset\lambda\,.
\endcases
$$
Here $H(\mu)=\prod\big(\mu_i-j+\theta(\mu_j'-i)+1\big)$, the
product runs over all boxes of the Young diagram $\mu$, the
length of the $i$-th row is denoted by $\mu_i$, and $\mu_j'$ is
the length of the $j$-th column of $\mu$. There is also an
explicit combinatorial formula for $\theta$-shifted Jack
polynomials,
$$
P_\mu^*(x;\theta) = \sum_T \psi_T(\theta)
\prod_{b\in\mu} \left(x_{T(b)} - c_\theta(b)\right),
\tag 6.2
$$
similar to the formula (A.7) for ordinary Jack polynomials. In
this formula $c_\theta(b)=(j-1)-\theta(i-1)$ stands for the
$\theta$-content of the box $b=(i,j)$ on the crossing of the
$i$-th row and $j$-th column of $\mu$. Note that $T$ ranges in
(6.2) over the set of {\it reverse} tableaux of shape $\mu$. See
\cite{17} for details and for the proof of the formula.

It follows directly from the equations (A.7) and (6.2) that the
leading symmetric term of the $\theta$-shifted Jack polynomial
$P_\mu^*$ equals the ordinary Jack polynomial, $[P_\mu^*]=P_\mu$.

We denote by $\dim_\theta(\mu,\nu)$ the dimension function (2.2)
with respect to the Jack multiplicity function (5.2), and we use
the abbreviation
$\dim_\theta\nu\equiv\dim_\theta(\varnothing,\nu)$. A nice {\it
hook formula} is available for $\theta$-dimension of true Young
diagrams:
$$
\dim_\theta\nu = {|\nu|! \over H_\theta(\nu)} = |\nu|!
\prod_{(i,j)\in\nu} \big((\nu_i-j)+(\nu_j'-i)\theta+1\big)^{-1}.
$$
This formula can be deduced from the results of Stanley
(\cite{22}, Theorem 5.4) and Macdonald (\cite{13}, VI.10). To
this end we remark that, by definition of dimension function,
$$
p_1^n(\,\cdot\,) = \sum_{\lambda\in\Bbb{Y}_n}
\dim_\theta\lambda\, P_\lambda(\,\cdot\,;\theta),
$$
hence
$$
\dim_\theta\lambda =
\frac{(p_1^n,P_\lambda)}{(P_\lambda,P_\lambda)}.
$$
Consider the representation $P_\lambda=\sum_\mu c_\mu p_\mu$ in
the basis of power sum functions. Then the coefficient
$c=c_{(1^n)}$ of $p_1^n$ equals the leading term in $X$ in the
formula (10.20) of \cite{13}, VI.10. It follows that the
numerator in the last formula is
$c\,(p_1^n,p_1^n)=c\,n!\,\alpha^n$. One should also use the
identity (10.16) in \cite{13} for
$b_\lambda=(P_\lambda,P_\lambda)^{-1}$. For another proof see
\cite{20}, Section 5.

The importance of $\theta$-shifted Jack polynomials for the
evaluation of Martin kernel of the Jack graph $\Bbb{J}(\theta)$
roots in the following basic dimension formula.

\proclaim{Theorem 6.1 \rm (\cite{20}, (5.2))}
Let $\mu$ and $\nu$ be arbitrary Young diagrams and $m=|\mu|$,
$n=|\nu|$. Then the Martin kernel (2.3) of the Young graph with
Jack edge multiplicities (5.2) can be written in the form
$$
\frac{\dim_\theta(\mu,\nu)}{\dim_\theta\nu} =
\frac{P^*_\mu(\nu;\theta)}{n(n-1)\dots(n-m+1)}.
\tag 6.3
$$
\endproclaim

\noindent {\bf Remarks.}
In the important particular case of $\theta=1$ the formula (6.3)
was first established in \cite{20}, Section 8. The polynomials
$s_\mu^*(x)=P_\mu^*(x;1)$ are called {\it shifted Schur
polynomials}, they were defined in \cite{20} by a simple
determinantal formula
$$
s^*_\mu(x_1,\,\ldots,x_n) = \dfrac
{\det[(x_i + n - i \downharpoonright \mu_j + n - j)]}
{\det[(x_i + n - i \downharpoonright n - j)]},
\tag 6.4
$$
similar to the well-known Weyl formula for the ordinary Schur
functions. Here and below we use the symbol
$$
(x\downharpoonright{m})=x(x-1)\dots(x-m+1)
$$ 
for the falling factorial powers of $x$.

Since there is an explicit expression for $\dim_\theta\nu$,
formula (6.3) yields an expression for the number
$\dim_\theta(\mu,\nu)$. In case of $\theta=1$, we obtain (thanks
to the determinantal formula (6.4)) a simple algebraic expression
for the number $\dim_1(\mu,\nu)$ of standard Young tableaux of
the {\it skew} shape $\nu\setminus\mu$. "

\subhead 7. The asymptotics of $\theta$-shifted symmetric
polynomials \endsubhead
The second major ingredient in our proof of Theorem B is the
asymptotic formula (7.4) to be proved in the present Section.

Given a box $b=(i,j)$ the number
$$
c_\theta(b) = (j - 1) - \theta(i - 1)
\tag 7.1
$$
is referred to as {\it $\theta$-content} of the box $b$. A box
$b$ is said to be {\it positive} or {\it negative} according to
the sign of its $\theta$-content. More precisely, $b$ is positive
if $c_\theta(b)>0$, and negative if $c_\theta(b)\le0$.

We may consider a Young diagram $\nu=(\nu_1,\,\ldots,\nu_l)$ as a
collection of boxes,
$$
\nu \equiv \{b_{ij}:\;\; 1\le i\le l,\;\; 1\le j\le \nu_i\},
$$
and we split $\nu$ as a union of disjoint subsets of its positive
and negative boxes,
$$
\nu^+ = \{b \in\nu\mid c_\theta(b) > 0\}, \quad
\nu^- = \{b \in\nu\mid c_\theta(b)\le0\}.
$$
Denote by $r$ the number of rows in $\nu^+$, and by $s$ the
number of columns in $\nu^-$. Let
$$
a_1\ge a_2\ge \dots\ge a_r>0;\qquad
b_1\ge b_2\ge \dots\ge b_s>0
\tag 7.2
$$
denote the lengths of corresponding rows and columns. Clearly,
$$
\sum_{i=1}^r a_i + \sum_{j=1}^s b_j = n
$$
is the total number of boxes in the diagram $\nu$.

Given a Young diagram $\nu$, we associate with it a point
$$
\omega_\nu(\theta) = \left(
\frac{a_1}n\;,\dots,\frac{a_r}n\;,0,\dots\;;\;
\frac{b_1}n\;,\dots,\frac{b_s}n\;,0,\dots
\right)
\tag 7.3
$$
of the Thoma simplex $\Omega$.

Recall that to any symmetric polynomial $f$ there corresponds a
continuous function $\omega\mapsto f(\omega|\theta)$ on the Thoma
simplex $\Omega$ defined by the equation (5.4). We shall show
that the normalized value $F(\nu)/n^m$ of a $\theta$-shifted
symmetric polynomial $F\in\Lambda^\theta$ of degree $m$
at a Young diagram $\nu$ with $n$ boxes gets close to the value
$f(\omega_\nu(\theta)|\theta)$ of its leading symmetric term
$f=[F]$ at the point $\omega_\nu(\theta)\in\Omega$, as 
$n\to\infty$.

\proclaim{Theorem 7.1}
Denote by $\nu$ a Young diagram with $n$ boxes, and by
$F\in\Lambda^\theta$ a $\theta$-shifted symmetric polynomial of
degree $\deg{F}=m$ with the leading symmetric term
$f=[F]\in\Lambda$. Then
$$
\left|\frac{F(\nu)}{n^m} - f(\omega_\nu(\theta)|\theta)\right|
\le {C \over \sqrt{n}},
\tag 7.4
$$ 
where the constant $C$ depends on $F$ and $\theta$, but not on
$n$.
\endproclaim

\demo{Proof}
It suffices to prove (7.4) for the polynomials
$F=\widetilde{p}_m$, where
$$
\widetilde p_m(x) = \sum_{i=1}^\infty
\Big((x_i - \theta(i-1) \downharpoonright m) -
(-\theta(i-1) \downharpoonright m)\Big)
$$
is a variant of $\theta$-shifted power sum polynomial, similar to
that of (6.1). Note that the leading symmetric term of
$\widetilde p_m$ equals $[\widetilde p_m]=p_m$, the conventional
power sum symmetric function. The polynomials $p_m$ generate the
algebra $\Lambda$, hence every element $F\in\Lambda^\theta$ is a
polynomial in the variables $\widetilde p_1,\widetilde
p_2,\ldots$. Since both functions $\nu\mapsto{F}(\nu)/n^m$,
$\nu\mapsto{f}(\omega_\nu(\theta);\theta)$ are bounded on
$\Bbb{Y}$, the Theorem will follow.

Recall (see Section A.1) that $(x\downharpoonright{m})$ denotes
the descending factorial power of $x$. The basic idea of the
proof is to evaluate a sum
$$
S_m = m\; \sum_{b\in\nu} (c_\theta(b) \downharpoonright m-1)
$$
in two different ways (a similar approach was employed in
\cite{8}).

First of all we note that $S_m$ coincides with the value of
$\widetilde p_m$ at the diagram $\nu$. Indeed, since
$(c+1\downharpoonright{m})-(c\downharpoonright{m})=
m(c\downharpoonright{m-1})$, we get an elementary summation
formula
$$
m \sum_{c=c_{\min}}^{c_{\max}-1} (c \downharpoonright m-1) =
(c_{\max} \downharpoonright m) - (c_{\min} \downharpoonright m).
\tag 7.5
$$
Splitting the diagram $\nu$ into its rows, and applying (7.5) to
each row, we see that
$$
S_m = \sum_{i\ge1} m \sum_{b\in\nu_i}
(c_\theta(b) \downharpoonright m-1) =
\sum_{i\ge1} \left( (\nu_i-\theta(i-1) \downharpoonright m) -
(-\theta(i-1) \downharpoonright m) \right) =
\widetilde p_m(\nu).
$$

In our second calculation of the sum $S_m$ we find partial sums
along the rows of the positive part $\nu^+$, and along the
columns of the negative part $\nu^-$ of the diagram $\nu$.

\proclaim{Lemma 7.1}
For every $m=1,2,\ldots$,
$$
\left|{m \over n^m} \sum_{b\in\nu^+}
(c_\theta(b) \downharpoonright m-1) -
\sum_{i=1}^r \Big({a_i \over n}\Big)^m \right|
\le {C_1 \over \sqrt{n}},
$$
where the constant $C_1$ does not depend on $n$.
\endproclaim

\demo{Proof}
Denote by $c_i$ the $\theta$-content (7.1) of the leftmost box in
the $i$-th row $\nu_i^+$ of the shape $\nu^+$. By the summation
formula (7.5),
$$
S_m^+(\nu)\equiv m \sum_{b\in\nu_i^+}
(c_\theta(b) \downharpoonright m-1) =
(a_i + c_i \downharpoonright m) - (c_i \downharpoonright m).
\tag 7.6
$$
An important feature of the partition $\nu=\nu^+\bigcup\nu^-$ is
that the $\theta$-contents of ``diagonal'' boxes are uniformly
bounded, $-1<c_i\le1$. As a result, one can derive from the
equation (7.6) an estimate
$$
\left| m \sum_{b\in\nu_i^+} (c_\theta(b) \downharpoonright m-1) -
a_i^m \right| \le \text{\tenrm const } a_i^{m-1},
\tag 7.7
$$
where {\it const} depends on $m$, but not on $n$ and $i$.

Remark that the Young diagram $\nu$ always contains a rectangle
on the crossing of its first $r$ rows and $s$ columns. Since
$r\approx\theta s$, it follows that
$$
r \le \text{\tenrm const }\, \sqrt{n}, \qquad
s \le \text{\tenrm const }\, \sqrt{n},
\tag 7.8
$$
where {\it const} depends on $\theta$ only.

Now divide both sides of (7.7) by $n^m$ and take into account the
inequalities $a_i\le n$. We obtain an estimate
$$
\left|{m \over n^m} \sum_{b\in\nu_i^+}
(c_\theta(b) \downharpoonright m-1) -
\Big({a_i \over n}\Big)^m \right|
\le {\text{\tenrm const } \over n}
$$
uniform in $i$, and the Lemma follows from the inequalities (7.8).
\qed\enddemo

Let us now deal with the column sum
$$
S_m^-(\nu) =
m\, \sum_{b\in\nu^-} (c_\theta(b) \downharpoonright m-1).
$$
We reduce its evaluation to that of the row sum
$$
S_m^+(\mu) =
m\, \sum_{b\in\mu^+} (c_\theta(b) \downharpoonright m-1)
$$
for the transposed Young diagram $\mu=\nu'$. Remark that the
positive part $\mu^+$ of $\mu$ is now taken with respect to the
$1/\theta$-content. The transposition map $(i,j)\mapsto(j,i)$
provides a bijection $b\mapsto{b}'$ between the boxes in $\nu^-$
and those in $\mu^+$ (up to a minor asymmetry in the definitions
of $\mu^+$ and $\nu^-$ which is not essential for our purposes).
Clearly,
$$
c_\theta(b) = (j-1) - \theta(i-1) =
(-\theta) \big((i-1)-\theta^{-1}(j-1)\big) =
(-\theta)\, c_{1/\theta}(b').
$$

\proclaim{Lemma 7.2}
Let $\mu=\nu'$ denote the transposed diagram of the Young diagram
$\nu$. Then
$$
{m \over n^m}
\left|
\sum_{b\in\nu^-} (c_\theta(b) \downharpoonright m-1) -
(-\theta)^{m-1}
\sum_{b\in\mu^+} (c_{1/\theta}(b') \downharpoonright m-1)
\right|
\le {C_2 \over \sqrt{n}},
$$
where the constant $C_2$ depends on $m$ and $\theta$, but not on
$n$.
\endproclaim

\demo{Proof}
Note that $(-\theta c\downharpoonright m-1)$ and
$(-\theta)^{m-1}(c\downharpoonright m-1)$ are both polynomials of
degree $m-1$ in the variable $c$, with one and the same leading
term $(-\theta c)^{m-1}$. It follows that
$$
\left| \sum_{b\in\nu_j^-} \Big(
(c_\theta(b) \downharpoonright m-1) -
(-\theta)^{m-1} (c_{1/\theta}(b') \downharpoonright m-1)
\Big)\right| \le \text{\tenrm const } b_j^{m-1},
$$
where {\it const} does not depend on $n$. Using inequalities
$b_j\le n$, we get
$$
{m \over n^m}
\left| \sum_{b\in\nu_j^-} \Big(
(c_\theta(b) \downharpoonright m-1) -
(-\theta)^{m-1} (c_{1/\theta}(b') \downharpoonright m-1)
\Big)\right| \le {\text{\tenrm const } \over n},
$$
and the Lemma follows from the estimates (7.8).
\qed\enddemo

By Lemma 7.1, we obtain formulae
$$
{m \over n^m} \sum_{b\in\nu^+}
(c_\theta(b) \downharpoonright m-1) =
\sum_{i=1}^r \Big({a_i \over n}\Big)^m +
O\Big({1\over\sqrt{n}}\Big)
$$
and
$$
{m \over n^m} \sum_{b\in\mu^+}
(c_{1/\theta}(b) \downharpoonright m-1) =
\sum_{j=1}^s \Big({b_j \over n}\Big)^m +
O\Big({1\over\sqrt{n}}\Big),
$$
so that the Lemma 7.2 and the equation
$S_m=S_m^+(\nu)+S_m^-(\nu)=\widetilde p_m(\nu)$ imply
$$
{\widetilde p_m(\nu) \over n^m} =
\sum_{i=1}^r \Big({a_i \over n}\Big)^m +
(-\theta)^{m-1}\, \sum_{j=1}^s \Big({b_j \over n}\Big)^m
+ O\left({1 \over \sqrt{n}}\right).
$$
which completes the proof of the Theorem.
\qed\enddemo

In the course of the proof of Theorem 7.1 we have used the
notation $\omega_\nu(\theta)$ for a point of the Thoma simplex
introduced by equation (7.3). It depends on the parameter
$\theta$ and differs from the point $\omega_\nu=\omega_\nu(1)$
introduced in a similar way by equation (3.3). Let us now remark
that if $\theta\ne0$ and $n=|\nu|$ goes to infinity, the elements
$\omega_\nu(\theta)$ and $\omega_\nu$ are asymptotically
equivalent. In fact, for each fixed $k=1,2,\ldots$ the
corresponding coordinates only differ by a constant (depending on
$k$).

\proclaim{Corollary 7.1}
Given a sequence $\nu^{(n)}\in\Bbb{Y}_n$ of Young diagrams,
$n=1,2,\ldots$, assume that their images (3.3) in the Thoma
simplex converge to a point $\omega\in\Omega$. Then, in the
notations of Theorem 7.1,
$$
\lim_{n\to\infty} {F(\nu^{(n)}) \over n^m} = f(\omega|\theta).
\tag 7.9
$$
\endproclaim

\subhead 8. Regular sequences of Young diagrams \endsubhead
In this Section we identify the Martin compactification of the
Jack graph $\Bbb{J}(\theta)$, $\theta>0$, with the space
$$
\widetilde{\Omega} = \big(\{0\} \times \Omega \big) \cup
\widetilde{\Bbb{Y}}\;\; \subset\;\; [0,1]\times\Omega
\tag 8.1
$$
(considered as a sort of geometric compactification). Recall that
the discrete subset $\widetilde{\Bbb{Y}}\subset[0,1]\times\Omega$
was introduced in Section 3 by the equations (3.3) and (3.4).
Both spaces $\widetilde{\Bbb{Y}}$ and $\widetilde{\Omega}$ were
defined with no reference to $\theta$.

\proclaim{Theorem 8.1}
Let $\nu^{(n)}$, $n=1,2,\ldots$ be a sequence of Young diagrams.
Then the following two conditions are equivalent:

(i)\; the sequence $(|\nu^{(n)}|^{-1},\nu^{(n)})$ converges in
$\widetilde{\Omega}$, as $n\to\infty$, to a point
$\widetilde{\omega}\in\widetilde{\Omega}$;

(ii)\; for each Young diagram $\mu$ there exists the limit
$K(\mu,\widetilde{\omega})=
\underset{n\to\infty}\to\lim K(\mu,\nu^{(n)})$.

\noindent
The Martin kernel of the Jack graph $\Bbb{J}(\theta)$ is provided
by the $\theta$-extended versions of Jack symmetric polynomials,
$$
K(\mu,\omega) = P_\mu(\omega|\theta).
\tag 8.2
$$
Here $P_\mu(\omega|\theta)=
\varphi_{\alpha,\beta}\big(P_\mu(\,\cdot\,;\theta)\big)$ denotes
the value of $\theta$-extended Jack polynomial indexed by $\mu$
at the point $\omega=(\alpha;\beta)$ of Thoma simplex, see
Definition (5.3).
\endproclaim

\demo{Proof}
Let us write $N(n)=|\nu^{(n)}|$ for the number of boxes in a
diagram $\nu^{(n)}$. Since the set $\widetilde{\Bbb{Y}}$ is
discrete, there will be no loss of generality in assuming that
$\lim N(n)=\infty$. Assume also that the condition (i) holds, so
that the sequence $\omega_{\nu^{(n)}}$ converges in the Thoma
simplex $\Omega$ to a point $\omega=(\alpha;\beta)$. By
definition of convergence in $\Omega$ this means that
$$
\lim_{n\to\infty} {(\nu^{(n)})_k \over N(n)} = \alpha_k;\qquad
\lim_{n\to\infty} {(\nu^{(n)})'_k \over N(n)} = \beta_k;\qquad
k = 1,2,\ldots\;.
\tag 8.3
$$

By the equation (6.3) we know that
$$
K(\mu,\nu^{(n)}) = {P_\mu^*(\nu^{(n)};\theta) \over N(n)^m} +
O\Big({1 \over N(n)}\Big),
$$
for all $\mu\in\Bbb{Y}_m$, $m=1,2,\ldots$. Applying Corollary 7.1
with $F(\nu)=P^*_\mu(\nu,\theta)$, we derive that
$$
\lim_{n\to\infty} K(\mu,\nu^{(n)}) =
P_\mu(\omega|\theta)
$$
for every $\mu\in\Bbb{Y}$. 

In the opposite direction, assume that the condition (ii) holds.
Since the space $\Omega$ is compact, one can choose a subsequence
of the sequence $\nu^{(n)}$, converging to a point
$\omega\in\Omega$. But the Jack polynomials form a linear basis
in the algebra $\Lambda$, and the $\theta$-extended symmetric
functions corresponding to polynomials in $\Lambda$ are uniformly
dense in the space $C(\Omega)$ of continuous functions on the
Thoma simplex (by Lemma 5.3). Hence, all partial limits of the
sequence $\nu^{(n)}$ are equal, and the condition (i) holds, too.
\qed\enddemo

The sequences $\nu^{(n)}\in\Bbb{Y}_n$ of Young diagrams subject
to the condition (ii) of Theorem 8.1 are called {\it regular
sequences}.

\proclaim{Corollary 8.1}
Assume that a regular sequence $\nu^{(n)}\in\Bbb{Y}_n$ converges to
a point $\omega\in\Omega$. Then the limit
$$
\lim_{n\to\infty}
{\dim_\theta(\mu,\nu^{(n)}) \over \dim_\theta \nu^{(n)}} =
P_\mu(\omega|\theta)
\tag 8.4
$$
exists for all $\mu\in\Bbb{Y}$.
\endproclaim

\proclaim{Corollary 8.2}
For every point $\omega\in\Omega$ the function
$\mu\mapsto{P}_\mu(\omega|\theta)$ is non-negative, normalized,
and harmonic with respect to the Jack multiplicities (5.2).
\endproclaim

\demo{Proof}
The harmonicity condition (2.1) coincides with the Pieri
formula (5.1) for Jack symmetric polynomials. The positivity
follows from Corollary 8.1, since the dimension functions are
non-negative.
\qed\enddemo

The remaining part of the proof of Theorems A, B is quite
standard and general; it does not depend on the specific features
of Jack graphs.

\proclaim{Lemma 8.1}
Let $\phi\in\Cal{H}$ be a normalized non-negative harmonic
function. Then
$$
M_n(\lambda) = \dim_\theta\lambda\;\phi(\lambda), \qquad
\lambda \in \Bbb Y_n
\tag 8.5
$$
is a probability distribution.
\endproclaim

\demo{Proof}
By the harmonicity condition (2.1), 
$$
M_{n-1}(\lambda) = \sum_{\nu:\lambda\nearrow\nu}
\frac{\dim_\theta\lambda\; \varkappa_\theta(\lambda,\nu)}
{\dim_\theta\nu}\; M_{n}(\nu),
$$
for every diagram $\lambda\in\Bbb{Y}_{n-1}$. By the definition
(2.2) of the dimension function,
$$
\sum_{\lambda:\lambda\nearrow\nu}
\frac{\dim_\theta\lambda\; \varkappa_\theta(\lambda,\nu)}
{\dim_\theta\nu} = 1,
$$
which implies
$$
\sum_{\nu\vdash n} M_n(\nu) =
\sum_{\lambda\vdash (n-1)} M_{n-1}(\lambda),
$$
and hence $\underset{\nu\vdash{n}}\to\sum M_n(\nu)=1$ for all
$n=1,2,\ldots$.
\qed\enddemo

\demo{Proof of Theorem B}
By the Corollary 8.2, a function
$$
\varphi(\mu) = \int_\Omega K(\mu,\omega)\, M(d\omega); \qquad
\mu \in \Bbb Y,
\tag 8.6
$$
is an element of $\Cal{H}$, for every probability distribution
$M$ on the Thoma simplex $\Omega$. By Lemma 5.3 such an integral
representation is unique. It only remains to find a measure $M$
representing a given harmonic function $\varphi$.

We show that for every function $\varphi\in\Cal{H}$ there exists
a weak limit $M=\lim\widetilde{M}_n$, where $\widetilde{M}_n$
denotes the image of the discrete probability distribution (8.5)
with respect to the embedding (3.3) of the set $\Bbb{Y}_n$ into
the Thoma simplex $\Omega$. The measure $M$ will provide the
representation (8.6) of the function $\varphi$.

In fact, it follows from (2.1) and (8.5) that
$$
\varphi(\mu) = \sum_{\nu\in\Bbb{Y}_n}
{\dim_\theta(\mu,\nu) \over \dim_\theta\nu} M(\nu),
$$
for every fixed $\mu\in\Bbb{Y}_m$ and $n>m$. Using equations
(6.3) and (7.4) we derive that
$$
\varphi(\mu) = \int_\Omega
P_\mu(\omega|\theta)\, \widetilde{M}_n(d\omega) + O(1/\sqrt{n}).
\tag 8.7
$$
Let us choose such a subsequence that the measures
$\widetilde{M}_n$ converge to a limiting measure $M$. Then
according to the formula (8.7)
$$
\varphi(\mu) = \int_\Omega P_\mu(\omega|\theta)\, M(d\omega),
$$
and the measure $M$ represents the harmonic function $\varphi$.
Note that the integral representation is unique, hence the limit
$\lim\widetilde{M}_n=M$ does exist.
\qed\enddemo

\demo{Proof of Theorem A}
It is clear that every point $(\alpha,\beta)\in\Omega$ determines
a positive evaluation homomorphism
$\varphi_{\alpha,\beta}(f)=f(\alpha;\beta)$, $f\in\Lambda$.

By Theorem B, every positive homomorphism
$\varphi:\Lambda\to\Bbb{R}$ (normalized by the condition
$\varphi(p_1)=1$) can be written in the form
$\varphi(f)=\int{f}\,dM$, for appropriate probability
distribution $M$ on the Thoma simplex $\Omega$. We have to show
that $M$ is a $\delta$-measure. Assuming this is not true, there
exists a symmetric function $f\in\Lambda$ with a non-trivial
distribution with respect to the measure $M$. But the variance of
$F$ vanishes,
$$
\text{\rm var}\,(f) =
\int_\Omega f^2\,dM - \Big(\int_\Omega f\,dM\Big)^2 =
\varphi(f^2) - \varphi^2(f) = 0,
$$
so that the assumption leads to a contradiction. (The argument is
a part of the proof of Theorem 6 in \cite{10}.) This implies
Theorem A.
\qed\enddemo

\head Appendix \endhead

We recall in this Section the basic background notation and
terminology related to integer partitions and symmetric
functions, see \cite{13} for a detailed exposition.
 
\subhead A.1. Young diagrams \endsubhead
Every decomposition
$\lambda:\;n=\lambda_1+\lambda_2+\ldots+\lambda_l$, where
$\lambda_1\ge\lambda_2\ge\ldots\ge\lambda_l$ are positive
integers, is called an {\it integer partition} of {\it length}
$l=l(\lambda)$. We identify integer partitions with corresponding
Young diagrams, and we denote by $n=|\nu|$ the number of boxes in
$\lambda$. The set of Young diagrams with $n$ boxes is denoted as
$\Bbb{Y}_n$, and $\Bbb{Y}=\bigcup\Bbb{Y}_n$ is the set of all
Young diagrams. We write $m_j=m_j(\lambda)$ for the number of
terms in $\lambda$ equal to $j$, and we set
$z_\lambda=\prod_{j\ge1}j^{m_j}\,m_j!$. The transposed
(conjugate) diagram of $\nu$ is denoted by $\nu'$; by definition,
$\nu_j'=\#\{i:\nu_i\ge{j}\}$.

The Young diagrams in $\Bbb{Y}$ are naturally ordered by
inclusion, and we shall write $\lambda\nearrow\nu$ to indicate
that the diagram $\nu$ covers $\lambda$ with respect to this
order. Another useful ordering of Young diagrams, denoted
$\lambda\ge\nu$, is called {\it dominance order}. By definition,
$\lambda\ge\nu$ iff $|\lambda|=|\nu|$ and
$\lambda_1+\ldots+\lambda_k\ge\nu_1+\ldots+\nu_k$ for all
$k\ge1$.

We use the symbol $(x\downharpoonright{m})=x(x-1)\dots(x-m+1)$
for the descending factorial powers of a variable $x$.
 
\subhead A.2. Symmetric functions \endsubhead
The $\Bbb{R}$-algebra of all symmetric polynomials in the
variables $x=(x_1,x_2,\ldots)$ is denoted by $\Lambda$. To each
monomial $x_1^{r_1}x_2^{r_2}\ldots$ of degree $n$ there
corresponds a partition $n=(\lambda_1,\lambda_2,\ldots)$, the
rearrangement of the exponents $r_1,r_2,\ldots$ in the decreasing
order. The monomial symmetric function $m_\lambda(x)$ is defined
as the sum of all monomials with a fixed partition $\lambda$. In
particular, the functions
$$
p_m(x) = m_{(n)}(x) = \sum_{j=1}^\infty x_j^m
$$
are called the {\it power sum symmetric polynomials}, and we set
$p_\lambda(x)=\prod_{j\ge1}p_j^{m_j(\lambda)}$, for every integer
partition $\lambda$. Each one of the families $\{m_\lambda\}$,
$\{p_\lambda\}$ forms a linear basis in the algebra $\Lambda$.

\subhead A.3. Extended symmetric functions \endsubhead
The {\it extended Schur functions} $s_\lambda(\alpha,\beta)$ can
be formally defined by the Jacoby -- Trudi determinant
$$
s_{(\lambda_1,\,\ldots,\lambda_m)} =
\left|\matrix
h_{\lambda_1}&h_{\lambda_1+1}&h_{\lambda_1+2}&\ldots&h_{\lambda_1+m-1}\\
h_{\lambda_2-1}&h_{\lambda_2}&h_{\lambda_2+1}&\ldots&h_{\lambda_2+m-2}\\
h_{\lambda_3-2}&h_{\lambda_3-1}&h_{\lambda_3}&\ldots&h_{\lambda_3+m-3}\\
   \ldots      &   \ldots      & \ldots      &\ldots& \ldots \\
h_{\lambda_m-m+1}&h_{\lambda_m-m+2}&h_{\lambda_m-m+3}&\ldots&h_{\lambda_m}
\endmatrix\right|,
\tag A.1
$$
where the {\it extended complete homogeneous symmetric functions}
$h_n=h_n(\alpha,\beta)$ arise as the coefficients of the
generating series
$$
e^{z\gamma}\prod_{j=1}^\infty {1 + z\beta_j \over 1 - z\alpha_j}
= 1 + \sum_{k=1}^\infty h_n(\alpha,\beta)\; z^n.
\tag A.2
$$

More generally, an {\it extended version} of a symmetric
polynomial may be defined as follows. Recall that every symmetric
function $f=f(\alpha)$ is a polynomial in the power sum symmetric
functions $p_m(\alpha)=\sum \alpha_j^m$, $m=1,2,\ldots$. In order
to obtain the extended version $f(\alpha;\beta)$ of the function
$f$, one only has to replace each occurrence of $p_1(\alpha)$ in
the corresponding polynomial by $1$, and each occurrence of
$p_m(\alpha)$, $m\ge2$, by the {\it extended power sum
polynomials}
$$
p_m(\alpha;\beta) = \sum_{j=1}^\infty \alpha_j^m +
(-1)^{m-1} \sum_{j=1}^\infty \beta_j^m.
\tag A.3
$$

\subhead A.4. Jack polynomials \endsubhead
Let us recall a few basic facts related to Jack symmetric
polynomials $P_\mu(x;\theta)$. Consider a scalar product on the
algebra $\Lambda$ defined by the equation
$$
(p_\mu,p_\lambda)_\theta = \delta_{\mu,\lambda}\,
z_\lambda\,\theta^{-l(\lambda)}
\tag A.4
$$
(note that we use the parameter $\theta$ inverse to that
employed by Macdonald in \cite{13}).

The polynomials $P_\mu(x;\theta)$ can be characterized by the
following properties:
$$
\gather
\big(P_\lambda(\,\cdot\,;\theta),P_\nu(\,\cdot\,;\theta)\big)_\theta=0
\qquad\text{\tenrm unless } \lambda=\nu;
\tag A.5 \\
P_\lambda(x;\theta) = m_\lambda(x) +
\sum_{\mu:\mu<\lambda} c_{\lambda\mu}(\theta)\, m_\mu(x),
\tag A.6
\endgather
$$
where $c_{\lambda\mu}(\theta)$ are rational functions of
$\theta$. In particular, the polynomials
$P_{(1^m)}(x;\theta)=m_{(1^m)}(x)$ do not depend on $\theta$.

There is a number of equivalent definitions of Jack polynomials;
for instance,
$$
P_\mu(x;\theta) = \sum_T \psi_T(\theta) \prod_{b\in\mu} x_{T(b)},
\tag A.7
$$
where $\psi_T(\theta)$ is an explicit rational expression in
$\theta$, see \cite{13}, VI.10.12.

In the limit $\theta\to 0$ the polynomials $P_\mu(x;0)$
degenerate to monomial symmetric functions $m_\mu(x)$, see
\cite{13}, p.380. The polynomials $s_\mu(x)=P_\mu(x;1)$ coincide
with the Schur symmetric functions.


\bigskip
\bigskip
{\smc S.~Kerov}: Steklov Mathematical Institute, St.~Petersburg
Branch (POMI), \newline Fontanka 27, 191011 St.~Petersburg,
Russia. E-mail address: {\tt kerov\@pdmi.ras.ru}
\newline\indent
{\smc A.~Okounkov}: Dept. of Mathematics, University of Chicago,
5734 South University Av., Chicago, IL 60637-1546, U.S.A. E-mail
address: {\tt okounkov\@math.uchicago.edu}
\newline\indent
{\smc G.~Olshanski}: Institute for Problems of Information
Transmission, Bolshoy \newline Karetny 19, 101447 Moscow GSP-4,
Russia. E-mail address: {\tt olsh\@ippi.ac.msk.su}


\Refs\widestnumber\key{20}

\ref\key 1 
\by M.~Aissen, I.~J.~Schoenberg, A.~M.~Whitney
\paper On the generating function of totally positive sequences I
\jour J. Analyse Math. 
\vol 2
\yr 1952
\pages 93--103
\endref

\ref\key 2 
\by J.~L.~Doob
\paper Discrete potential theory and boundaries
\jour J. of Math. and Mech.
\vol 8
\yr 1959
\pages 433-458
\endref

\ref\key 3 
\by A.~Edrei
\paper On the generating function of totally positive sequences II
\jour J. Analyse Math. 
\vol 2
\yr 1952
\pages 104--109
\endref

\ref\key 4 
\by V.~N.~Ivanov
\paper Dimensions of skew shifted Young diagrams and projective
characters of the infinite symmetric group
\inbook Representation Theory, Dynamical Systems, Combinatorial
and Algorithmical Methods, II
\bookinfo Zapiski Nauchnykh Seminarov POMI, Vol. 240
\yr 1997
\pages 116--136
\lang Russian
\transl\nofrills English transl. to appear in J.~Math.~Sciences
\endref

\ref\key 5 
\by S.~V.~Kerov
\paper Combinatorial examples in the theory of AF-algebras
\inbook Differential geometry, Lie groups and mechanics X
\bookinfo Zapiski Nauchnykh Seminarov LOMI, Vol. 172
\yr 1989
\pages 55--67
\lang Russian
\transl\nofrills English transl. in  J.~Soviet Math. {\bf 59} 
(1992), No. 5, pp.~1063--1071.
\endref

\ref\key 6 
\by S.~V.~Kerov
\paper Generalized Hall-Littlewood symmetric functions and orthogonal
polynomials
\inbook Representation Theory and Dynamical Systems
\ed A.~M.~Vershik
\bookinfo Advances in Soviet Math., Vol. 9
\publ Amer. Math. Soc.
\publaddr Providence, R.I.
\yr 1992
\pages 67--94
\endref

\ref\key 7 
\by S.~V.~Kerov
\paper The Boundary of Young Lattice and Random Young Tableaux
\jour DIMACS Series in Discrete Mathematics and
Theoretical Computer Science
\vol 24 
\yr 1996
\pages 133--158
\endref

\ref\key 8 
\by S.~Kerov  and G.~Olshanski
\paper Polynomial functions on the set of Young diagrams
\jour Comptes Rend. Acad. Sci. Paris, Ser. I
\vol 319
\yr 1994
\pages 121--126
\endref

\ref\key 9 
\by S.~Kerov, A.~Vershik
\paper The Grothendieck Group of the Infinite Symmetric Group and
Symmetric Functions with the Elements of the $K_0$-functor theory
of AF-algebras
\inbook Representation of Lie groups and related topics
\bookinfo Adv. Stud. Contemp. Math.
\vol 7
\eds A.~M.~Vershik and D.~P.~Zhelobenko
\publ Gordon and Breach
\yr 1990
\pages 36--114
\endref

\ref\key 10
\by S.~V.~Kerov, A.~M.~Vershik
\paper Characters, factor-representations and  $K$-functor of the
infinite symmetric group
\inbook Proc. of the International Conference on Operator Algebras and
Group Representations
\vol 1
\yr 1980
\pages 23--32
\endref

\ref\key 11 
\by J.~F.~C.~Kingman
\paper Random partitions in population genetics
\jour Proc. Roy. Soc. London A.
\vol 361
\yr 1978
\pages 1--20
\endref

\ref\key 12 
\by F.~Knop and S.~Sahi
\paper Difference equations and symmetric polynomials
defined by their zeros
\jour Internat.\ Math.\ Res.\ Notices 
\yr 1996 \issue 10 \pages 473--486
\endref 

\ref\key 13 
\by I.~G.~Macdonald
\book Symmetric functions and Hall polynomials
\bookinfo 2nd edition
\publ Oxford University Press
\yr 1995
\endref

\ref\key 14 
\by A.~I.~Markushevich
\book Theory of Analytic Functions
\lang Russian
\publaddr Moscow -- Leningrad
\yr 1950
\endref

\ref\key 15 
\by M.~L.~Nazarov
\paper Projective representations of the infinite symmetric group
\inbook Representation Theory and Dynamical Systems
\ed A.~M.~Vershik
\bookinfo Advances in Soviet Math., Vol. 9
\publ Amer. Math. Soc.
\publaddr Providence, R.I.
\pages 115--130
\yr 1992
\endref

\ref\key 16 
\by A.~Yu.~Okounkov
\paper Thoma's theorem and representations of infinite bisymmetric
group 
\jour Funct. Analysis and its Appl.
\vol 28
\yr 1994
\pages 101--107
\endref

\ref\key 17 
\by A.~Yu.~Okounkov
\paper (Shifted) Macdonald polynomials:
$q$-integral representation and combinatorial formula
\jour q-alg/9605013
\endref

\ref\key 18 
\by A.~Yu.~Okounkov
\paper On representations of the infinite symmetric group
\inbook Representation Theory, Dynamical Systems, Combinatorial
and Algorithmical Methods, II
\bookinfo Zapiski Nauchnykh Seminarov POMI, Vol. 240
\yr 1997
\pages 167--229
\lang Russian
\transl\nofrills English transl. to appear in J.~Math.~Sciences
\endref

\ref\key 19 
\by A.~Okounkov and G.~Olshanski
\paper Shifted Schur functions
\jour Algebra i Analiz
\vol 9
\yr 1997
\pages No.~2
\lang Russian
\transl\nofrills English version to appear in St.~Petersburg Math. J. 
{\bf 9} (1998), No.~2;  q-alg/9605042
\endref

\ref\key 20 
\bysame
\paper Shifted Jack polynomials, binomial formula,
and applications
\jour Math. Research Lett. 
\vol 4
\yr 1997
\pages 69--78 (see also q-alg/9608020)
\endref

\ref\key 21 
\by G.~I.~Olshanski
\paper Unitary representations of $(G,K)$-pairs connected with the
infinite symmetric group $S(\infty)$
\jour Leningrad Math. J. 
\vol 1
\yr 1990
\pages 983--1014
\endref

\ref\key 22 
\by R.~P.~Stanley
\paper Some Combinatorial Properties of Jack Symmetric Functions
\jour Adv. Math.
\vol 77
\yr 1989
\pages 76-115
\endref

\ref\key 23 
\by E.~Thoma
\paper Die unzerlegbaren, positive-definiten Klassenfunktionen
der abz\"ahlbar unendlichen, symmetrischen Gruppe
\jour Math.~Zeitschr.
\vol 85
\yr 1964
\pages 40-61
\endref

\ref\key 24 
\by A.~M.~Vershik, S.~V.~Kerov
\paper Asymptotic character theory of the symmetric group
\jour Funct. Analysis and its Appl.
\vol 15
\yr 1981
\pages 246--255
\endref

\endRefs
\enddocument
\end